# Online Learning Based Risk-Averse Stochastic MPC of Constrained Linear Uncertain Systems

Chao Ning, Fengqi You*

College of Engineering, Cornell University, Ithaca, New York 14853, USA
*Corresponding author. Phone: (607) 255-1162; Fax: (607) 255-9166; E-mail: fengqi.you@cornell.edu

**Abstract**

This paper investigates the problem of designing data-driven stochastic Model Predictive Control (MPC) for linear time-invariant systems under additive stochastic disturbance, whose probability distribution is unknown but can be partially inferred from data. We propose a novel online learning based risk-averse stochastic MPC framework in which Conditional Value-at-Risk (CVaR) constraints on system states are required to hold for a family of distributions called an ambiguity set. The ambiguity set is constructed from disturbance data by leveraging a Dirichlet process mixture model that is self-adaptive to the underlying data structure and complexity. Specifically, the structural property of multimodality is exploited, so that the first- and second-order moment information of each mixture component is incorporated into the ambiguity set. A novel constraint tightening strategy is then developed based on an equivalent reformulation of distributionally robust CVaR constraints over the proposed ambiguity set. As more data are gathered during the runtime of the controller, the ambiguity set is updated online using real-time disturbance data, which enables the risk-averse stochastic MPC to cope with time-varying disturbance distributions. The online variational inference algorithm employed does not require all collected data be learned from scratch, and therefore the proposed MPC is endowed with the guaranteed computational complexity of online learning. The guarantees on recursive feasibility and closed-loop stability of the proposed MPC are established via a safe update scheme. Numerical examples are used to illustrate the effectiveness and advantages of the proposed MPC.

## 1. Introduction

Over the past few decades, Model Predictive Control (MPC) has established itself as a modern control strategy with theoretical grounding and a wide variety of applications (Mayne, 2014; Mayne et al., 2000; Morari & Lee, 1999; Qin & Badgwell, 2003) because it addresses multivariate dynamic systems subject to control input and state constraints. Implemented in a receding horizon fashion, MPC solves a finite-horizon optimal control problem at each sampling instant and only performs the first control action. This procedure is repeated at the next instant with a new measurement update. However, the presence of uncertainty could inflict severe performance degradation or even loss of feasibility on conventional MPC if the uncertainty is not explicitly accounted for (Rawlings & Mayne, 2009).

Motivated by this fact, much research has been done on designing MPC that accounts for the uncertainty of predictions. Robust MPC strategies, in which disturbances are modeled using a bounded and deterministic set, aim to satisfy the hard constraints of states and control inputs for all possible uncertainty realizations (Bemporad & Morari, 1999; Chisci et al., 2001; Langson et al., 2004; Mayne et al., 2005). Designing robust MPC is necessary and efficient when state and input constraints need to be immunized against uncertainty. By exploiting the probabilistic nature of uncertainty in the controller design, stochastic MPC methods are capable of tolerating constraint violations in a systematic way (Mesbah, 2016). Additionally, stochastic MPC can increase the region of attraction by means of chance constraints, thus allowing for a systematic trade-off between constraint satisfaction and control performance (Farina et al., 2016).

Due to this attractive feature, stochastic MPC has stimulated considerable research interest from the control community (Farina, et al., 2016; Mayne, 2016; Mesbah, 2016). The existing literature in stochastic MPC can be typically grouped into two main categories depending on whether the knowledge of uncertainty distribution is perfectly known or not. The first category subsumes the research that focuses on tightening or adaptive tightening for chance constraints via an explicit use of uncertainty distributions (Cannon et al., 2011; Chatterjee et al., 2011; Chatterjee & Lygeros, 2015; Dai et al., 2015; Korda et al., 2014; Kouvaritakis et

---

The material in this paper was not presented at any conference.

al., 2010; Lorenzen et al., 2017a; Muñoz-Carpintero et al., 2018). By considering stability and feasibility requirements separately, a novel constraint tightening method was proposed for stochastic MPC based on a disturbance distribution known to significantly increase a feasible region (Lorenzen, et al., 2017a). Conventional stochastic MPC strategies rely heavily on the assumption that the probability distribution of disturbance is known *a priori*. However, such perfect knowledge of distribution is rarely available in practice. Instead, only partial information regarding disturbance distribution can be inferred from data. In such a data-driven setting, chance constraints in these stochastic MPC frameworks are no longer ensured because there always exists a gap between the underlying true distribution and an estimated one. For the case with an unknown disturbance distribution, the probabilistic constraints can be reformulated using a scenario-based approach (Calafiore & Fagiano, 2013; Lorenzen et al., 2017b; Schildbach et al., 2014), or by means of the Chebyshev-Cantelli inequality (Farina et al., 2015; Farina & Scattolini, 2016; Paulson & Mesbah, 2017). To speed up computation, a novel offline sampling-based approach, along with its sample complexity, was developed for discrete-time linear systems subject to multiplicative stochastic disturbance (Lorenzen, et al., 2017b). In addition, the Chebyshev-Cantelli inequality-based method guarantees constraint satisfaction for any distributions sharing the same mean-covariance information, which is similar to the emerging paradigm called distributionally robust control (Singh et al., 2019; Sopasakis et al., 2019; Van Parys et al., 2016; Yang, 2018). In the paradigm of distributionally robust control, an ambiguity set is a set comprising all possible probability distributions characterized by certain statistical properties of stochastic disturbances (Van Parys, et al., 2016). Distributionally robust control methodologies have been proposed for linear systems with probabilistic constraints assuming the first- and second-order moment information of uncertainty distribution. However, only global mean and covariance information is utilized, and the existing ambiguity set is not updated in an online manner based on newly collected disturbance data.

With the ever-increasing availability of data in control systems, there is a growing trend to leverage the information embedded within data to improve control performance. The remarkable progress in machine learning and big data analytics leads to a broad range of opportunities to integrate data-driven systems with model-based control systems (Shang & You, 2019), enabled by the dramatic growth in computing power. Recently, learning-based MPC has attracted increasing attention from the control community (Aswani et al., 2013; Koller et al., 2018; Rosolia & Borrelli, 2018; Rosolia et al., 2018). One such method leveraged statistical identification tools to build a data-driven system model for enhanced performance, while utilizing an approximate model with uncertainty bounds for constraint tightening in robust tube MPC (Aswani, et al., 2013; Limon et al., 2017). Along the same research direction, a learning-based robust MPC was developed to integrate control design with offline system model learning via a set membership identification (Terzi et al., 2019). For a single-input single-output system under stochastic uncertainty, an adaptive dual MPC was designed (Heirung et al., 2017), where control played the role of a probe to learn about the system model (Filatov & Unbehauen, 2000). To learn system nonlinearities from data, several MPC strategies leverage Gaussian process regression which provides system dynamic models as well as residual bounds (Hewing et al., 2017; Soloperto et al., 2018). Learning-based MPC is well suited for the exploitation of data values to enhance control performance (Rosolia & Borrelli, 2018), and as such makes it a practical and appealing control tool in the era of big data (Wu et al., 2019). A comprehensive review of learning-based MPC was provided, including the research on learning the system dynamics for MPC, learning the controller design, and safe learning (Hewing et al., 2020). However, most of existing learning-based MPC methods focus on learning system models, which essentially learn a function by means of regression techniques along with their error bounds. Although some of these learning-based MPC methods allow for online or adaptive model learning, most of them focus on the robust control framework. Few studies have organically integrated online learning with stochastic MPC. Therefore, from both a theoretical and practical standpoint, a systematic investigation of online learning based stochastic MPC with theoretical guarantees on recursive feasibility and closed-loop stability is needed.

To fill this research gap, there are several computational and theoretical challenges that need to be addressed. The first research challenge is how to devise a data-driven ambiguity set of disturbance distributions in a nonparametric manner such that it is adaptive towards data complexity automatically. In the context of online learning, one cannot easily pin down the complexity of a disturbance model at the beginning, as more data stream in over the runtime of MPC and data complexity can grow over time. Another key research challenge is how to develop a framework that organically integrates online learning with MPC for intelligent control. In particular, with more and more data collected, the online learning method needs to be scalable with sample size in terms of both memory and computational time. The third challenge lies in the development of a computationally tractable constraint tightening method which hedges against distributional ambiguity in stochastic MPC while leveraging structural properties of disturbance distributions. This calls for the theoretical extension in distributionally robust optimization with the nonparametric ambiguity set, because there are no theoretical results available for direct use. The fourth key research challenge is how to guarantee recursive feasibility and stability of online learning based stochastic MPC. This challenge arises from the integration between online learning with stochastic MPC. Specifically, the online update of disturbance distribution information could jeopardize the property of recursive feasibility.

This paper proposes an online learning based risk-averse stochastic MPC framework for linear time-invariant systems affected by additive disturbance. Instead of assuming perfect knowledge of a disturbance distribution, we consider a more realistic setting where the distribution can be partially inferred from data. To immunize the control strategy against distributional ambiguity, Conditional Value-at-Risk (CVaR) constraints are required to be satisfied for all candidate distributions in an ambiguity set. As opposed to the conventional ambiguity set using global mean-covariance information, a nonparametric data-driven ambiguity set is constructed by taking advantage of a Dirichlet Process Mixture Model (DPMM) (Gomes et al., 2008). Specifically, we devise the ambiguity set based on a structural property, namely multimodality, along with local first- and second-order moment information of each mixture component, the number of which is automatically derived from disturbance data. During the runtime of the controller, real-time disturbance data are exploited to update the disturbance distribution based on an online variational inference algorithm. The developed distribution-learning-while-control scheme alternates between learning the ambiguity set from real-time disturbance data and controlling the system with updated uncertainty information. Afterwards, we propose a constraint tightening technique for the resulting distributionally robust CVaR constraints over the DPMM-based ambiguity set. Specifically, the constraints are equivalently reformulated as Linear Matrix Inequality (LMI) constraints, which are amenable to efficient computation. Additionally, we introduce a safe online update scheme for the ambiguity set such that recursive feasibility and closed-loop stability are ensured. Numerical simulation and comparison studies show that the proposed MPC method enjoys less-conservative control performance compared with the conventional distributionally robust control that uses global mean and covariance information of disturbance. Additionally, thanks to the online learning scheme, the proposed MPC is advantageous in terms of a low constraint violation percentage under time-varying disturbance distributions.

The major contributions of this paper are summarized as follows:
- A novel online Bayesian learning based risk-averse stochastic MPC framework that improves control performance based on real-time data;
- An online data-driven approach with the DPMM to devise ambiguity sets that are self-adaptive to the underlying structure and complexity of real-time disturbance data;
- A novel constraint tightening method for risk-averse stochastic MPC, in which data-driven CVaR constraints over the DPMM-based ambiguity set are equivalently reformulated as computationally tractable LMIs;
- Theoretical guarantees on recursive feasibility and stability of the proposed MPC with the introduction of a novel safe update scheme for ambiguity sets.

One important contribution of this work is a novel and organic integration of online learning with stochastic MPC under possibly time-varying disturbance distributions. Moreover, a novel nonparametric ambiguity set based on the DPMM is first developed in this paper. To the best of our knowledge, this learning based ambiguity set is the first of its kind that is self-adaptive to the underlying disturbance data structure and complexity. Note that the objective of the constraint tightening method is to facilitate the solution process of the proposed MPC. To the best of our knowledge, no existing publication addresses the theoretical guarantees on recursive feasibility and stability of risk-averse stochastic MPC in the face of time-varying disturbance distributions that can be inferred from data. Therefore, the establishment of recursive feasibility and stability with the safe update scheme is another novel contribution.

The remainder of this paper is structured as follows. The problem setup and preliminaries are presented in Section 2. An online learning based stochastic MPC framework is then proposed in Section 3. The theoretical properties, including recursive feasibility and closed-loop stability are derived in Section 4. Numerical examples are presented in Section 5 to demonstrate the effectiveness and advantages of the proposed approach, followed by a concluding remark.

*Notation:* The notation used in this paper is standard. For sets $\mathbb{A}$ and $\mathbb{B}$, $\mathbb{A} \oplus \mathbb{B} = \{a+b | a \in \mathbb{A}, b \in \mathbb{B}\}$ and $\mathbb{A} \ominus \mathbb{B} = \{a | a \oplus \mathbb{B} \subseteq \mathbb{A}\}$ denote the Minkowski set addition and Pontryagin set difference, respectively. $[A]_j$ represents the $j$-th row of a matrix $A$, and $[a]_j$ denotes the $j$-th entry of a vector $a$. For matrices $A$ and $B$, their Frobenius inner product is represented as $A \bullet B$. The concatenated vector of $c_{0|k}, \ldots, c_{N-1|k}$ is denoted by $\mathbf{c}_{N|k}$. For integers $p$ and $q$, the notation $[p, q]$ represents all the integer numbers between $p$ and $q$. For a real number $\alpha$, $(\alpha)^+ = \max(\alpha, 0)$.

## 2. Problem setup and preliminaries

In this section, we describe the problem setup, including system dynamics, distributionally robust CVaR constraints on system states, hard constraints on control inputs, an objective function, as well as some basics on risk-averse stochastic MPC.

Consider the following linear, time-invariant uncertain system with additive stochastic disturbance.

$$x_{k+1} = Ax_k + Bu_k + w_k \tag{1}$$

where $x_k \in \mathbb{R}^n$ is a system state, $u_k \in \mathbb{R}^m$ denotes a control input, and $w_k \in \mathbb{R}^n$ represents an additive disturbance. Let $\mathbb{X} = \{x \in \mathbb{R}^n | Hx \leq h\}$ and $\mathbb{U} = \{u \in \mathbb{R}^m | Gu \leq g\}$ be the polytopic constraints on the system states and the control inputs, both of which contain the origin in the interior.

We make the following assumptions regarding the system and the additive disturbance.

**Assumption 1.** The measurement of system state $x_k$ is available at time $k$.

This is a common assumption. At each sampling time $t+1$, the realization of disturbance at previous time instant can be obtained by $w_t = x_{t+1} - Ax_t - Bu_t$. Therefore, we have access to real-time disturbance data for online learning which is used to refine the knowledge of uncertainty.

**Assumption 2.** The matrix pair $(A, B)$ is controllable.

**Assumption 3.** The additive disturbance $w$ has a bounded and convex support set $\mathbb{W} = \{w | Ew \leq f\}$, which contains the origin in the interior.

Given the state measurement at time $k$, the predicted model in MPC is given by,

$$x_{l+1|k} = Ax_{l|k} + Bu_{l|k} + w_{l|k}, \quad x_{0|k} = x_k \quad (2)$$

where $x_{l|k}$ and $u_{l|k}$ represent the $l$-step ahead system state and control input predicted at time $k$, respectively. Disturbance $w_{l|k}$ denotes the $l$-step ahead stochastic disturbance with distributional properties available at time $k$.

By leveraging the probability distribution of disturbance $w$, stochastic MPC allows constraint violation via the use of chance constraints. Following the literature (Korda et al., 2011; Lorenzen, et al., 2017a; Paulson et al., 2018), we consider the chance constraints of states, as in (3). The advantage of such "one-step-ahead" chance constraints is that it facilitates recursive feasibility analysis and ensures closed-loop performance.

$$\mathbb{P}\{[H]_i x_{l+1|k} \leq [h]_i | x_{l|k}\} \geq 1 - [\varepsilon]_i, \quad i \in [1, p] \quad (3)$$

where $H \in \mathbb{R}^{p \times n}$, $h \in \mathbb{R}^p$, and parameter $[\varepsilon]_i$ is a pre-specified risk level for $i$-th constraint on the predicted system state. Note that (3) becomes hard constraints when the probability mass of all disturbances is strictly greater than zero and $[\varepsilon]_i$ is set to be zero.

Chance constraints ensure the constraints are satisfied with a probability of at least $1-[\varepsilon]_i$, relying on the assumption that the disturbance distribution is perfectly known. Stochastic MPC can enlarge the feasible region of the corresponding finite-horizon optimal control problem via tunable risk levels in chance constraints, thus improving objective performance. By assuming a disturbance distribution, tightening parameters can be computed offline by the inverse of the cumulative distribution of disturbance. In practice, however, such precise knowledge of the disturbance distribution is rarely available, and only partial information can be inferred from historical as well as real-time disturbance data. Due to the finite amount of disturbance data, the assumed disturbance probability $\mathbb{P}$ could deviate from the underlying true distribution. Consequently, the actual constraint violation resulting from conventional stochastic MPC could become worse than the pre-specified one. Additionally, the chance constraints *per se* (in the form of (3)) focus on the frequency of constraint violation and fail to account for violation magnitude without penalizing the magnitude in a cost function as implemented in soft constrains (Samuelson & Yang, 2018; Van Parys, et al., 2016). Therefore, we introduce the definition of CVaR and distributionally robust CVaR constraints as follows.

**Definition 1** (*Conditional Value-at-Risk*). For a given measureable loss $L: \mathbb{R}^n \to \mathbb{R}$, a random vector $\tilde{\xi}$ on $\mathbb{R}^n$ with probability $\mathbb{P}$, and tolerance level $\epsilon \in (0,1)$, the CVaR of random loss function $L$ at level $\epsilon$ with respect to probability $\mathbb{P}$ is defined below.

$$\mathbb{P}\text{-CVaR}_\epsilon \left(L(\tilde{\xi})\right) = \inf_{\beta \in \mathbb{R}} \left\{\beta + \frac{1}{\epsilon} \mathbb{E}_\mathbb{P}\left\{\left(L(\tilde{\xi}) - \beta\right)^+\right\}\right\} \quad (4)$$

where $\mathbb{E}_\mathbb{P}$ denotes the expectation with respect to probability $\mathbb{P}$. The CVaR can be interpreted as the conditional expectation of loss $L$ above the $(1-\epsilon)$ quantile of the probability distribution of $L$.

To address the issue of conventional chance constraints, we propose to using a distributionally robust CVaR version of constraints $\mathbb{P}\{[H]_i x_{l+1|k} \leq [h]_i | x_{l|k}\} \geq 1 - [\varepsilon]_i$, which is provided as follows.

$$\sup_{\mathbb{Q} \in \mathcal{Q}^{(k)}} \mathbb{Q}\text{-CVaR}_{[\varepsilon]_i} \left([H]_i x_{l+1|k} - [h]_i\right) \leq 0 \quad (5)$$

where $\mathcal{Q}^{(k)}$ is defined as an ambiguity set constructed based on disturbance data information available up to time $k$, and $\mathbb{Q}$ denotes the conditional probability given $x_{l|k}$. Unlike the chance constraints, the distributionally robust CVaR constraint $\sup_{\mathbb{Q} \in \mathcal{Q}^{(k)}} \mathbb{Q}\text{-CVaR}_{[\varepsilon]_i} \left([H]_i x_{l+1|k} - [h]_i\right) \leq 0$ not only hedges against distributional ambiguity, but also penalizes severe constraint violations that could be detrimental to control systems.

Meanwhile, hard constraints are imposed on control inputs due to the physical limitation of actuators, as shown below.

$$Gu_{l|k} \leq g \quad (6)$$

where $G \in \mathbb{R}^{q \times m}$ and $g \in \mathbb{R}^q$.

Following a common practice in stochastic MPC, we split the predicted state into its nominal part and stochastic error part as follows:

$$x_{l|k} = z_{l|k} + e_{l|k} \quad (7)$$

where $z_{l|k}$ and $e_{l|k}$ denote the nominal part and stochastic error part of the predicted state $x_{l|k}$, respectively.

An MPC strategy with a prediction horizon $N$ is considered. By employing an error feedback, the predicted input for the uncertain system can be represented by

$$u_{l|k} = Ke_{l|k} + v_{l|k} \quad (8)$$

where $K$ is a stabilizing feedback gain, and $v_{l|k}$ represents the predicted control input for the nominal system corresponding to $z_{l|k}$. The nominal control input $v_{l|k}$ can be formulated as follows.

$$v_{l|k} = Kz_{l|k} + c_{l|k} \tag{9}$$

where $c_{l|k} \in \mathbb{R}^m$, $l = 0,...,N-1$ are decision variables in the receding horizon optimal control problem, and $c_{l|k} = 0$, $l \geq N$.

With predictive control laws $u_{l|k} = Ke_{l|k} + v_{l|k}$ and $v_{l|k} = Kz_{l|k} + c_{l|k}$, the system dynamics for the nominal state and stochastic error are given below.

$$z_{l+1|k} = Az_{l|k} + Bv_{l|k}, \quad z_{0|k} = x_k \tag{10}$$

$$e_{l+1|k} = \Phi e_{l|k} + w_{l|k}, \quad e_{0|k} = 0 \tag{11}$$

where $\Phi = A + BK$ is Schur stable with the stabilizing feedback gain $K$.

We consider the cost for the nominal system in this work. Specifically, the control objective is to minimize the infinite horizon cost at sampling time $k$, as given below.

$$J_\infty = \sum_{l=0}^{\infty} \left( z_{l|k}^T Q z_{l|k} + v_{l|k}^T R v_{l|k} \right) \tag{12}$$

where $Q \in \mathbb{R}^{n \times n}$, $Q \succ 0$, $R \in \mathbb{R}^{m \times m}$, and $R \succ 0$.

Furthermore, the following assumption of detectability is made such that there exists an LQR solution.

**Assumption 4.** The matrix pair $(A, Q^{1/2})$ is detectable.

Suppose the feedback gain matrix $K$ is chosen to be LQR optimal, we can further rewrite (12) as (13). The detailed derivation is provided in Appendix A.

$$J_\infty = \sum_{l=0}^{N-1} c_{l|k}^T \left( R + B^T PB \right) c_{l|k} + x_k^T P x_k \tag{13}$$

where $P$ denotes the solution of the Lyapunov equation $\Phi^T P \Phi + Q + K^T RK = P$. Note that the second term $x_k^T P x_k$ is a constant. A quadratic finite-horizon cost is given by

$$J_N\left(\mathbf{c}_{N|k}\right) = \sum_{l=0}^{N-1} c_{l|k}^T \left( R + B^T PB \right) c_{l|k} \tag{14}$$

where $\mathbf{c}_{N|k} = \left( c_{0|k}^T, ..., c_{N-1|k}^T \right)^T$ represents the decision vector at time $k$ with a horizon of length $N$.

## 3. Online learning based risk-averse stochastic MPC

In this section, we propose a novel risk-averse stochastic MPC framework based on online Bayesian learning. We first develop a data-driven approach to construct an ambiguity set for the stochastic disturbance based on the DPMM. Then, an efficient constraint tightening method for the CVaR constraints on system states over the ambiguity set is developed for the synthesis of a stochastic predictive controller. Finally, based on an online safe update scheme, the predictive control algorithm that organically integrates online learning with risk-averse stochastic MPC is described.

### 3.1. Online Bayesian learning for streaming disturbance measurement data

To automatically decipher the structural property of the disturbance distribution, we employ a nonparametric Bayesian model, called the DPMM, which is briefly described as follows.

A Dirichlet Process (DP) constitutes a fundamental building block for the DPMM that relies on mixtures to characterize data distribution. The DP is technically a probability distribution over distributions. Suppose a random distribution $G$ follows a DP parameterized by a concentration parameter $\alpha$ and a base measure $G_0$ over space $\Theta_0$, denoted as $G \sim DP(\alpha, G_0)$. For any fixed partitions $(A_1, ..., A_r)$ of $\Theta_0$, we have the following

$$\left(G(A_1),...,G(A_r)\right) \sim Dir\left(\alpha G_0(A_1),...,\alpha G_0(A_r)\right) \tag{15}$$

where $Dir$ represents a Dirichlet distribution.

Following a stick-breaking procedure (Sethuraman, 1994), a random draw from the DP can be expressed by $G = \sum_{k=1}^{\infty} \pi_k \delta(\tilde{\varphi}_k)$, where $\pi_k = \bar{\beta}_k \prod_{j=1}^{k-1} (1 - \bar{\beta}_j)$ is the weight, $\tilde{\varphi}_k$ is sampled from $G_0$, and $\delta(\tilde{\varphi}_k)$ denotes the Dirac delta function at $\tilde{\varphi}_k$. The parameter $\bar{\beta}_k$ represents the proportion being broken from the remaining stick, and follows a Beta distribution, denoted as $\bar{\beta}_k \sim Beta(1, \alpha)$.

The Bayesian nonparametric model, i.e. the DPMM, employs $\tilde{\varphi}_k$ as parameters of some data distribution. Based on the DP, we summarize the basic form of the DPMM as follows (Blei & Jordan, 2006; Campbell & How, 2015):

$$\begin{aligned} &\{\pi_k, \tilde{\varphi}_k\}_{k=1}^{\infty} \sim DP(\alpha, G_0) \\ &l_n \sim Mult(\pi) \\ &o_n \sim F(\tilde{\varphi}_{l_n}) \end{aligned} \tag{16}$$

where $Mult$ denotes a multinomial distribution, $l_n$ is the label indicating the component or cluster of observation $o_n$, $n$ is an index ranging from 1 to $N_d$, and data $o_1, ..., o_{Nd}$ are distributed according to distribution $F$. Based on the stick-breaking procedure, $G$ is discrete with probability one. Such discreteness further induces the clustering of data.

Due to its computational efficiency, variational inference has become a method of choice for approximating the conditional distribution of latent variables in the DPMM given observed data (Blei & Jordan, 2006). In variational inference, the problem of computing the posterior distribution is formulated as an optimization problem, which can be solved using a coordinate ascent method. Following the literature (Campbell & How, 2015), we use the mixtures of Gaussian distributions in this paper. Therefore, we can choose $\tilde{\varphi}_n \sim NW\left(\theta_k, \lambda_k, \omega_k, \Psi_k^{-1}\right)$, where $\tilde{\varphi}_k = \left(\eta_k, \tilde{H}_k\right)$ includes mean vector and precision matrix, and $NW$ represents a normal Wishart distribution. A variational distribution is used to approximate the true posterior in terms of Kullback–Leibler divergence (Blei & Jordan, 2006; Kurihara et al., 2007).

For the online learning setting, suppose that the real-time data $w_t$ is collected for the control system. To learn from the streaming data, we employ an online variational inference algorithm in this work (Gomes, et al., 2008). This algorithm features faster computation and bounded memory requirement for each round of learning. It is well-suited to learning the distribution online from real-time data over the runtime of MPC. The algorithm iterates between a model building phase and a compression phase. In the model building phase, clump constraint $C_s$ is the set of indices satisfying that $\forall i, j \in C_s$, and disturbance data $w_i$ and $w_j$ are generated from the same mixture component. Disturbance data within the same clump are summarized via the average sufficient statistics, which encapsulate all the information needed for the purpose of inference (Gomes, et al., 2008). The new disturbance data at time $t$ are used to update the inference results, then they are discarded to reduce the memory overhead. As a result, the algorithm is attractive in terms of both bounded memory requirement and fast computation. By introducing the compression phase, the algorithm not only is computationally efficient, but also requires bounded memory space. The clump constraints are determined in the compression phase in a top-down recursive fashion. Specifically, the computational burden at each model update using new disturbance data does not grow with the processed disturbance data amount. For more details on this online learning algorithm for the DPMM, we refer the readers to (Gomes, et al., 2008).

Based on the online variational inference results available up to time $k$, we propose a data-driven ambiguity set, as given in **Definition 2**, by leveraging both multimodality and moment information of each mixture component.

**Definition 2.** The ambiguity set based on the DPMM, denoted as $\mathcal{D}$, is defined as follows.

$$\mathcal{D} = \sum_{j=1}^{m^{(k)}} \gamma_j^{(k)} \mathcal{D}_j \left( \mathbb{W}, \mu_j^{(k)}, \Sigma_j^{(k)} \right) \tag{17}$$

where $\mathbb{W} = \{w | Ew \leq f\}$ is the support set of the additive disturbance under **Assumption 3**, $m^{(k)}$ denotes the number of mixture components, the mixing weight $\gamma_j^{(k)}$ indicates the occurring probability of each mixture component, $\mathcal{D}_j$ represents a basic ambiguity set, $\mu_j^{(k)}$ and $\Sigma_j^{(k)}$, respectively, denote the mean and covariance estimates for the $j$-th mixture component obtained from the online learning algorithm. Note that the support set plays a key role in ensuring recursive feasibility by means of a terminal set.

The data-driven ambiguity set for the stochastic disturbance based on the DPMM is devised as a weighted Minkowski sum of several basic ambiguity sets, the number of which is automatically determined from disturbance data using the online variational inference algorithm. Each basic ambiguity set $\mathcal{D}_j$ is cast as follows.

$$\mathcal{D}_j \left( \mathbb{W}, \mu_j^{(k)}, \Sigma_j^{(k)} \right) = \left\{ \rho \in \mathcal{M}_+ \left| \begin{array}{l} \int_{\mathbb{W}} \rho(d\xi) = 1 \\ \int_{\mathbb{W}} \xi \cdot \rho(d\xi) = \mu_j^{(k)} \\ \int_{\mathbb{W}} \xi \xi^T \cdot \rho(d\xi) \leq \Sigma_j^{(k)} + \mu_j^{(k)} \left( \mu_j^{(k)} \right)^T \end{array} \right. \right\} \tag{18}$$

where $\mathcal{M}_+$ represents the set of positive Borel measures on $\mathbb{R}^n$, and $\rho$ is a positive measure.

There are several highlights of the proposed data-driven ambiguity set. First, the ambiguity set is devised in a nonparametric manner so that it automatically accommodates its complexity to the underlying structure and complexity of disturbance data. Second, each basic ambiguity set is devised using mean and covariance information, which endows the resulting stochastic MPC with enormous computational benefits. Third, the proposed ambiguity set leverages the fined-grained distribution information, namely local moment information. This feature implies that the resulting stochastic MPC enjoys less conservative control performance compared with a control method utilizing a conventional ambiguity set based on global moment information.

*3.2. A novel constraint tightening method based on distributionally robust CVaR constrained optimization*

In this section, we present a novel constraint tightening technique for distributionally robust CVaR constraints $\sup_{\mathbb{Q} \in \mathcal{Q}^{(k)}} \mathbb{Q}\text{-CVaR}_{[\varepsilon]_i} \left( [H]_i x_{l+1|k} - [h]_i \right) \leq 0$ over the data-driven ambiguity set $\mathcal{D}$, as well as constraint tightening for control input constraints $Gu_{l|k} \leq g$. The purpose of constraint tightening on system states is to obtain constraints for states of the nominal system, such that the constraints on states of the uncertain system are satisfied with a frequency of at least $1-[\varepsilon]_i$.

For the state constraints, the corresponding distributionally robust CVaR constraints can be equivalently reformulated as LMI constraints on the predicted nominal system states, as given in the following theorem.

**Theorem 1.** (Constraint tightening). The system satisfies the data-driven distributionally CVaR constraints if and only if the nominal system satisfies the tightened constraints $z_{l+1|k} \in \tilde{\mathbb{Z}}^{(k)} \ominus \left( \bigoplus_{i=1}^{l} \Phi^i \mathbb{W} \right)$, with $\tilde{\mathbb{Z}}^{(k)}$ given below.

$$\tilde{\mathbb{Z}}^{(k)} = \left\{ z \in \mathbb{R}^n \middle| Hz \leq h - \eta^{(k)} \right\} \tag{19}$$

where $[\eta^{(k)}]_i$ is the optimal objective value of the following optimization problem (20).

$$\min \eta$$

s.t. $[\varepsilon]_i \cdot \beta_i + \sum_{j=1}^m \gamma_j \left\{ t_{ij} + \mu_j^T \omega_{ij} + \left( \Sigma_j + \mu_j \mu_j^T \right) \bullet \Omega_{ij} \right\} \leq 0$

$$\begin{pmatrix} \Omega_{ij} & \frac{1}{2}(\omega_{ij} + E^T \varphi_{ij}) \\ \frac{1}{2}(\omega_{ij} + E^T \varphi_{ij})^T & t_{ij} - f^T \varphi_{ij} \end{pmatrix} \geq 0, \forall j \quad (20)$$

$$\begin{pmatrix} \Omega_{ij} & \frac{1}{2}(\omega_{ij} - [H]_i^T + E^T \phi_{ij}) \\ \frac{1}{2}(\omega_{ij} - [H]_i^T + E^T \phi_{ij})^T & t_{ij} + \beta_i + \eta - f^T \phi_{ij} \end{pmatrix} \geq 0, \forall j$$

$$\Omega_{ij} \geq 0, \varphi_{ij} \geq 0, \phi_{ij} \geq 0, \forall j$$

Note that we drop the index $k$ from $m^{(k)}$, $\gamma_j^{(k)}$, $\mu_j^{(k)}$ and $\Sigma_j^{(k)}$ for notational simplicity. The constraints in (20) are LMI constraints. The proof of **Theorem 1** is provided in Appendix B.

**Remark 1**. The support set is assumed to be a polyhedron in **Assumption 3**, because it facilitates the LMI reformulation in the proposed constraint tightening approach. Specifically, (B.8) and (B.9) with a polyhedral support set can be reformulated as computationally efficient LMI constraints in the proposed MPC. In principle, the assumption on the polyhedral support can be relaxed to a general compact convex set. In this case, those constraints still admit robust counterparts, but could be less computationally efficient than LMIs (Ben-Tal et al., 2015).

The optimization problem in (20) is a computationally tractable Semi-Definite Program (SDP), which can be solved efficiently using off-the-shelf optimization solvers, such as SeDuMi and MOSEK. The distribution information of disturbance is learned online through the online variational inference algorithm, and is further incorporated to (20) to perform constraint tightening, thus improving control performance in an online fashion.

For hard constraints on control inputs $Gu_{l|k} \leq g$, the corresponding constraint tightening result is obtained as follows, employing a tube-based strategy.

$$v_{l|k} \in \mathbb{V}_l = \mathbb{U} \ominus K \left( \bigoplus_{i=0}^{l-1} \Phi^i \mathbb{W} \right) \quad (21)$$

where $\mathbb{U} = \left\{ u \in \mathbb{R}^m \mid Gu \leq g \right\}$.

### 3.3. The proposed online learning based risk-averse stochastic MPC algorithm with a safe update scheme

In this section, we first introduce the finite horizon optimal control problem of the proposed MPC with a safe ambiguity set update scheme. Then, the overall description of the resulting online learning based risk-averse stochastic MPC algorithm is provided in detail.

In the online learning based risk-averse stochastic MPC paradigm, the finite horizon optimal control problem is solved online repeatedly. The optimal control problem needed to be solved at time $k$ is denoted as (**OL-SMPC**$_k$), given in (22)-(27).

$$\min_{\mathbf{c}_{N|k}} J_N \left( \mathbf{c}_{N|k} \right) \quad (22)$$

s.t. $z_{l+1|k} = Az_{l|k} + Bv_{l|k}, \quad z_{0|k} = x_k \quad (23)$

$v_{l|k} = Kz_{l|k} + c_{l|k} \quad (24)$

$z_{l+1|k} \in \mathbb{Z}_{l+1}^{(k)}, \quad l \in [0, N-1] \quad (25)$

$v_{l|k} \in \mathbb{V}_l, \quad l \in [0, N-1] \quad (26)$

$z_{N|k} \in \mathbb{Z}_f^{(k)} \quad (27)$

where sets $\mathbb{Z}_{l+1}^{(k)}$ and $\mathbb{Z}_{l+1}^{(k)}$ are defined using the safe update scheme described as follows.

In the developed update scheme, the condition is identified under which it remains safe to utilize the tightened constraints with an updated ambiguity set in the risk-averse stochastic MPC. Specifically, if a candidate solution satisfies the tightened constraints resulting from current results of the online variational inference, one can safely incorporate the newly-learned uncertainty information into the control problem; otherwise, one resorts to the tightened constraints from the previous sampling time.

To describe the condition checked by the safe update scheme, we need the definition of the candidate solution given as follows.

**Definition 3.** (Candidate solution) Given an optimal solution $\mathbf{c}_{N|k}^* = \left( c_{0|k}^*, c_{1|k}^*, \ldots, c_{N-1|k}^* \right)$ to the MPC problem at time $k$, the candidate solution at time instant $k+1$ is defined by

$$\tilde{\mathbf{c}}_{k+1} = \left( c_{1|k}^*, \ldots, c_{N-1|k}^*, 0 \right) \quad (28)$$

Note that it is a typical approach to employ the candidate solution as the shifted optimal input augmented by zero (Brunner et al., 2017; Chisci, et al., 2001). This work employs a dual mode prediction paradigm (Kouvaritakis & Cannon, 2016), namely Mode 1 corresponds to $v_{l|k} = Kz_{l|k} + c_{l|k}$, $l = 0, \ldots, N-1$, and Mode 2 corresponds to $v_{l|k} = Kz_{l|k}$, $l \geq N$. As a result, the terminal controller $v_{l|k} = Kz_{l|k}$ for the nominal system $z_{l+1|k} = Az_{l|k} + Bv_{l|k}$ is able to steer the nominal system state (in the terminal set) to the origin. The last term of zero in the candidate solution actually comes from the shifted solution at Mode 2 ($c_{N|k} = 0$). Therefore, it is not restrictive for the last term to be zero. The explicitly given candidate solution plays a critical role not only in establishing recursive feasibility, but also in proving closed-loop stability (Brunner, et al., 2017; Chisci, et al., 2001).

Based on the candidate solution, the safe update scheme checks the following condition:

$$\tilde{z}_{l|k+1} \in \hat{\mathbb{Z}}_l^{(k+1)} \tag{29}$$

$$\tilde{z}_{N|k+1} \in \hat{\mathbb{Z}}_f^{(k+1)} \tag{30}$$

where $\tilde{z}_{l|k+1}$ denotes the predicted state of the nominal system corresponding to the candidate solution $\tilde{\mathbf{c}}_{k+1}$, set $\hat{\mathbb{Z}}_l^{(k+1)} = \tilde{\mathbb{Z}}^{(k+1)} \ominus \left( \bigoplus_{i=1}^{l} \Phi^i \mathbb{W} \right)$, and set $\hat{\mathbb{Z}}_f^{(k+1)}$ is defined using terminal constraints.

To formally define the terminal set $\hat{\mathbb{Z}}_f^{(k+1)}$, we define the (maximal) robust positively invariant set below (Blanchini, 1999).

**Definition 4.** (Robust positively invariant set) A set $\Omega$ is a robust positively invariant set for system $x_{k+1} = f(x_k, w_k)$ and constraint set $(\mathbb{X}, \mathbb{W})$ if $\Omega \subseteq \mathbb{X}$ and $x_{k+1} \in \Omega$, $\forall w_k \in \mathbb{W}$ for every $x_k \in \Omega$.

**Definition 5.** (Maximal robust positively invariant set) A set $\Omega$ is a Maximal Robust Positively Invariant (MRPI) set for system $x_{k+1} = f(x_k, w_k)$ and constraint set $(\mathbb{X}, \mathbb{W})$ if $\Omega$ is a robust positively invariant set and contains all robust positively invariant sets.

**Remark 2.** The MRPI sets are computed using a standard approach based on the recursion of predecessor sets, a.k.a. backward reachable sets (Blanchini, 1999; Kolmanovsky & Gilbert, 1998).

**Definition 6.** (Terminal set) The terminal set $\hat{\mathbb{Z}}_f^{(k)}$ is defined as the MRPI set for the following system $z_{N|k+1} = \Phi z_{N|k} + \Phi^N w_k$ that satisfies the tightened constraints $z_{N|k} \in \hat{\mathbb{Z}}_N^{(k)}$ and $K z_{N|k} \in \mathbb{V}_N$.

Based on **Definition 6**, the terminal set should be the MRPI set which satisfies $\Phi \hat{\mathbb{Z}}_f^{(k)} \oplus \Phi^N \mathbb{W} \subseteq \hat{\mathbb{Z}}_f^{(k)}$ and $\forall z \in \hat{\mathbb{Z}}_f^{(k)} \Rightarrow z \in \hat{\mathbb{Z}}_N^{(k)}$, $Kz \in \mathbb{V}_N$.

Note that the terminal controller using the feedback gain $K$ respects state and input constraints, when operated in the terminal set.

It is safe to update the ambiguity set when an indicator called *flag* equals 1; otherwise, when *flag*=0, updating the ambiguity set could jeopardize the recursive feasibility of the proposed MPC. Therefore, the corresponding tightened constraints adopted in the online learning based stochastic MPC are given by

$$\mathbb{Z}_l^{(k+1)} = flag \cdot \hat{\mathbb{Z}}_l^{(k+1)} \oplus (1-flag) \cdot \mathbb{Z}_l^{(k)} \tag{31}$$

Similarly, we can safely update the terminal set in the following way:

$$\mathbb{Z}_f^{(k+1)} = flag \cdot \hat{\mathbb{Z}}_f^{(k+1)} \oplus (1-flag) \cdot \mathbb{Z}_f^{(k)} \tag{32}$$

The proposed MPC algorithm is detailed in **Fig. 1**. The online learning based risk-averse stochastic MPC algorithm can be roughly divided into two blocks: (i) offline computation of sets, the ambiguity set construction based on historical disturbance data, and (ii) online learning from real-time disturbance data and online optimization. At each time, the MPC strategy only implements the first control action. From **Fig. 1**, we can see that it alternates the online optimal control of the system, and the online learning for real-time disturbance data. If the candidate solution satisfies the condition $\tilde{z}_{l|k+1} \in \hat{\mathbb{Z}}_l^{(k+1)}$ and $\tilde{z}_{N|k+1} \in \hat{\mathbb{Z}}_f^{(k+1)}$, the newly learned uncertainty information is incorporated into the predictive control strategy to improve control performance over its runtime. In the MPC framework, a finite-horizon optimal control problem is solved at each time step, and only the control decisions at the first step are implemented as control input $u_k$. The corresponding first control action $u_{0|k}^* = K e_{0|k} + v_{0|k}^*$. As $e_{0|k} = 0$, we can further have $u_{0|k}^* = 0 + v_{0|k}^* = v_{0|k}^*$. Thus, $u_k = v_{0|k}^*$ as in Step 5 of the algorithm.

| **Algorithm** Online-learning based stochastic MPC algorithm |
|---|
| 1: **Offline:** Given the initial state $x_0$, construct an ambiguity set in (17) from historical data, and determine $\mathbb{Z}_l^{(0)}$ using (19)-(20), $\mathbb{V}_l$ in (21), and $\mathbb{Z}_f^{(0)}$ based on Definition 6. |
| 2: **Online:** |
| 3: **for** $k$=0,1…**do** |
| 4:     Solve the optimal control problem (**OL-SMPC**$_k$) in (22)-(27); |
| 5:     Apply control policy in (8) for $l$=0, i.e. $u_k = v_{0|k}^*$; |
| 6:     Measure the current system state $x_{k+1}$ and obtain $w_k$ using (1); |
| 7:     Run online learning for DPMM in (16) with real-time data $w_k$; |
| 8:     Update ambiguity set in (17), and obtain $\hat{\mathbb{Z}}_l^{(k+1)}$ using (19)-(20) and $\hat{\mathbb{Z}}_f^{(k+1)}$ based on Definition 6; |
| 9:     **if** $\tilde{z}_{l|k+1} \in \hat{\mathbb{Z}}_l^{(k+1)}$ in (29) and $\tilde{z}_{N|k+1} \in \hat{\mathbb{Z}}_f^{(k+1)}$ in (30) |
| 10:       *flag*=1; |
| 11:     **else** |
| 12:       *flag*=0; |
| 13:     **end** |
| 14:     Determine set $\mathbb{Z}_l^{(k+1)}$ and $\mathbb{Z}_f^{(k+1)}$ using (31)-(32); |
| 15: **end** |

**Fig. 1.** The pseudocode of the proposed online learning based risk-averse stochastic MPC algorithm.

**Remark 3.** The upper bound of the memory requirement of the online learning algorithm is $\left( \frac{n^2 + 3n}{2} + 1 \right) N_c + n N_s$, where $N_c$ is the number of clumps, $N_s$ is the number of singlets, and $n$ denotes the data dimension. Note that computational cost for storing the sufficient statistics for each clump is $\left( \frac{n^2 + 3n}{2} + 1 \right)$. Compared with a batch learning algorithm, the computational complexity of the adopted online learning algorithm is only $O(K_m(N_c + N_s + 1))$ during the model building phase (Gomes, et al., 2008), where $K_m$ denotes the maximum number of components.

## 4. The theoretical properties of the proposed online learning based risk-averse stochastic MPC

In this section, the properties of the proposed online learning based risk-averse stochastic MPC algorithm, namely recursive feasibility (formally defined in **Definition 7**) and closed-loop stability, are established.

**Definition 7.** (Recursive feasibility) If the finite-horizon optimization problem of MPC is initially feasible, it remains feasible for all the subsequent sampling instant.

As pointed out in Section 1, an important property of MPC, namely recursive feasibility, might be compromised by the time-varying probabilistic constraints with an online updated ambiguity set of disturbance distributions. To this end, the novel safe update scheme for the ambiguity set is developed for the MPC framework, along with the terminal set or terminal constraints (Blanchini, 1999).

By employing the safe ambiguity set update scheme developed in Section 3.3, the recursive feasibility and closed-loop stability are ensured despite that the disturbance distribution might be time-varying. Note that we assume disturbance support $\mathbb{W}$ is time-invariant for the ease of exposition, as indicated in **Assumption 3** where matrices $E$ and $f$ are not indexed by time $k$. For the time-varying support, the derivation and the conclusion of the proposed constraint tightening approach in **Theorem 1** are still valid. The only issue with the varying support is that the MRPI set could become empty when the support becomes sufficiently large. This could lead to further infeasibility issues.

A standard assumption underpinning tube based MPC on the terminal set is made as follows (Langson, et al., 2004; Mayne, et al., 2005).

**Assumption 5.** There exists a nonempty terminal set $\mathbb{Z}_f^*$ for the tightened constraints based on a worst-case scheme, i.e. $\hat{\mathbb{Z}}_f^* = \tilde{\mathbb{Z}}^* \ominus \left( \bigoplus_{i=1}^l \Phi^i \mathbb{W} \right)$ with $\tilde{\mathbb{Z}}^* = \{ z \in \mathbb{R}^n | Hz \le h - \eta_0 \}$ where $[\eta_0]_i = \max_{\xi \in \mathbb{W}} [H]_i \xi$.

**Proposition 1.** Under **Assumption 5**, the terminal MRPI set $\mathbb{Z}_f^{(k)}$ is always nonempty with the constraint tightening approach in **Theorem 1**.

The proof of **Proposition 1** is given in Appendix C.

**Remark 4.** Based on the proof in Appendix C, we can see that set $\mathbb{Z}_f^*$ is a robust positively invariant set, not necessarily the MRPI set, for the updated tightened constraints. Therefore, for the case with a very limited computational budget, the offline-computed set $\mathbb{Z}_f^*$ can serve as terminal constraints to guarantee recursive feasibility and stability without the need of re-computing the MRPI set at each time step.

Next, we prove the recursive feasibility of the proposed MPC, as given in the following theorem.

**Theorem 2.** Let $\mathbb{D}(x_k)$ denote the feasible region of the finite horizon optimal control problem (**OL-SMPC$_k$**) for state $x_k$. If $\mathbb{D}(x_0) \ne \emptyset$, then given **Assumptions 1-5**, we have $\mathbb{D}(x_k) \ne \emptyset, \forall k$.

The proof of **Theorem 2** is provided in Appendix D.

**Remark 5.** The distribution of the stochastic disturbance can be arbitrarily time-varying, so the feasibility of the candidate solution based on the adaptive constraints (without using the safe update scheme) cannot hold universally. To the best of our knowledge, the proposed safe update scheme represents the first attempt to successfully address the recursive feasibility issue of stochastic MPC subject to time-varying disturbance distributions.

Before proving the stability of the proposed MPC, we provide the definition of minimal robust positively invariant set as follows.

**Definition 8.** The minimal robust positively invariant set $R_\infty$ is defined as follows.

$$R_\infty \doteq \lim_{l \to \infty} \bigoplus_{i=0}^l \Phi^i \mathbb{W} \quad (33)$$

The following theorem establishes the stability of the closed-loop system with the proposed MPC strategy.

**Theorem 3.** (Closed-loop stability). Given that $\mathbb{D}(x_0) \ne \emptyset$, the closed-loop system state asymptotically converges to a neighborhood of the original under the proposed online learning based risk-averse stochastic MPC.

The proof of **Theorem 3** is provided in Appendix E.

## 5. Numerical examples

In this section, we apply the proposed online learning based stochastic MPC to numerical examples to illustrate its effectiveness and advantages. We also implement the risk-averse stochastic MPC using global mean and covariance information in the ambiguity set and the risk-averse stochastic MPC without online learning, in addition to the proposed MPC approach for the purpose of comparison. The online learning algorithm and the MPC control strategies are implemented in MATLAB R2018a. We use the YALMIP toolbox in MATLAB R2018a (Lofberg, 2004). The GUROBI 8.0 solver is adopted to solve the finite-horizon optimal control problem, and SeDuMi 1.3 is employed to solve the constraint tightening problem (20). The computational experiments are performed on a computer with an Intel (R) Core (TM) i7-6700 CPU @ 3.40 GHz and 32 GB RAM. In practical applications, computational efficiency could be significantly improved with faster computing facilities and a better software environment. The related sets, including robust positively invariant sets, are obtained via Multi-Parametric Toolbox 3.0 (Herceg et al., 2013). The prediction horizon $N$ is set to be nine. Note that the distributionally robust control method (Van Parys et. al, 2016) considers a setting similar to ours, e.g. the risk-averse stochastic control setting where the disturbance distribution is only partially known. This is why we compare the proposed method with this existing distributionally robust control method. The comparison with other backoff adaptation methods, such as the sampling based method, is not apple-to-apple for the following reasons. First, although the sampling-based approach (Lorenzen, et al., 2017a, 2017b) also uses disturbance data, it assumes that the underlying disturbance distribution is fixed. However, this setting is different from the considered setting in our work, in which the disturbance distribution is only partially inferred from data and possibly time-varying. Second, the sampling approach merely provides an inner approximation for chance constraints, while the proposed approach develops an equivalent reformulation of distributionally robust CVaR chance constraints. Third, the uncertainty is required to be sampled offline to ensure recursive feasibility by suitably augmenting the constraints (Lorenzen, et al., 2017a, 2017b). By contrast, our proposed MPC leverages online learning from real-time disturbance data, and guarantees recursive feasibility and stability.

### 5.1. Example with a disturbance distribution having multimodality

In this section, we use a benchmark numerical example, control of a constrained sampled double integrator (Mayne, et al., 2005), to demonstrate the effectiveness of the proposed MPC. The system dynamics in this benchmark example is defined by

$$x_{k+1} = \begin{pmatrix} 1 & 1 \\ 0 & 1 \end{pmatrix} x_k + \begin{pmatrix} 0.5 \\ 1 \end{pmatrix} u_k + w_k \quad (34)$$

The state and control constraints in the risk-averse stochastic MPC is given as follows.

$$\sup_{\mathbb{P} \in \mathcal{D}^{(k)}} \mathbb{P}\text{-CVaR}_{0.2}\left([0\ 1]x_{l+1|k} - 2\right) \leq 0 \quad (35)$$

$$\mathbb{U} = \left\{ u \in \mathbb{R} \mid |u| \leq 5 \right\} \quad (36)$$

The initial condition of system states $x_0 = (-5, -2)^T$, and the support set of disturbance is $\mathbb{W} = \{w \mid \|w\|_\infty \leq 0.6\}$. In the numerical example, the matrix gain $K$ is determined to be the unconstrained LQR solution with matrix $Q = \begin{pmatrix} 1 & 0 \\ 0 & 1 \end{pmatrix}$ and $R=0.01$. Specifically, the matrix gain $K = [-0.6609\ -1.3261]$, so $\Phi = \begin{pmatrix} 0.6696 & 0.3370 \\ -0.6609 & -0.3261 \end{pmatrix}$. Set $R_\infty$ can be computed using an approximation method (Rakovic et al., 2005), in which an upper bound on the approximation error can be specified *a priori*.

To demonstrate the control performance and computational time, 100 closed-loop simulations are performed with a simulation horizon of 20 time steps. The closed-loop cost function is $J_{\text{cost}} = \sum_{k=1}^{T_s}\left(x_k^T Q x_k + u_k^T R u_k\right)$, and the simulation horizon length is $T_s=20$. Compared with the risk-averse stochastic MPC with a mean-covariance ambiguity set, the proposed MPC method exploits fine-grained uncertainty information, and is less conservative via reducing the closed-loop cost by an average of 9.77% over all simulation runs. To take a closer look at the computational time breakdown of the proposed MPC, we present the average computational time for online learning, constraint tightening through solving (20), online control via solving (**OL-SMPC**$_k$), and computing the MRPI set in **Fig. 2**. The average computational times are calculated over 2,000 (20 × 100) time steps. From the results, we can see that the average computational time for online learning is merely 0.344 s. It is computationally efficient to update pamater $\eta$ online, and it takes 0.121 s on average to solve problem in (20). The proposed MPC not only enables the incorporation of updated distribution information to improve control performance, but also features an acceptable computational cost. Because the example is a benchmark problem, the computational time comparison is representative of the computational performance for the proposed approach.

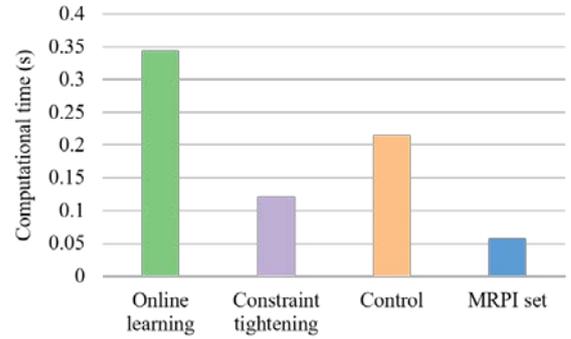

**Fig. 2.** The average computational times of the proposed online learning based risk-averse stochastic MPC method over 2,000 time steps.

### 5.2. Example with a time-varying disturbance distribution

To demonstrate the significance of online learning in the proposed MPC approach, we consider an example with a time-varying disturbance distribution. Specifically, the standard deviation used to generate historical disturbance data is 0.005, while the standard deviation increases to 0.3 for the real-time disturbance data generation. In this example, the system dynamics are provided in (34), the support set of disturbance is $\mathbb{W} = \{w \mid \|w\|_\infty \leq 0.1\}$, matrix $Q = \begin{pmatrix} 1 & 0 \\ 0 & 1 \end{pmatrix}$, $R=1$, $\mathbb{U} = \{u \in \mathbb{R} \mid |u| \leq 1\}$ and the distributionally robust CVaR constraint is given below.

$$\sup_{\mathbb{P} \in \mathcal{D}^{(k)}} \mathbb{P}\text{-CVaR}_{0.15}\left([0\ 1]x_{l+1|k} - 1.2\right) \leq 0 \quad (37)$$

To take a closer look at constraint violations under time-varying disturbance distributions, 100 simulations are performed with different realizations of disturbance sequences. **Fig. 3** shows a set of state trajectories using the proposed MPC for a simulation horizon of 20 time steps.

Note that the corresponding prediction horizon is $N=9$. For the proposed MPC strategy, the average constraint violation of the first nine time steps is 7.2%, even if the disturbance variation becomes significant because of the larger standard deviation. In constrast, the performance of risk-averse stochastic MPC using distributionally robust CVaR constraints without online learning deteriorates, and the corresponding average constraint violation increases to 22.3%. This constraint violation is higher than the prescribed tolerance of 15.0%. Notably, the risk-averse stochastic MPC without the online learning scheme implements the constraint tightening in (20) offline using historical disturbance data and its parameter $\eta$ remains constant over time.

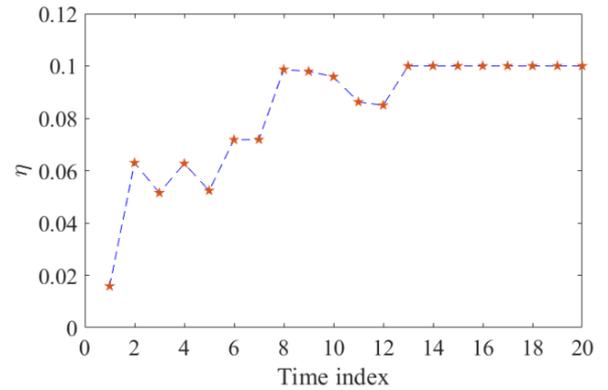

**Fig. 4.** The online adaption of constraint tightening parameters in the proposed MPC for a time-varying disturbance distribution in a simulation.

### 5.3. Discussions on MRPI set computation

In this section, we provide discussions on MRPI set computation and put forward a variant of the proposed MPC to address the issue of re-computing MPRI sets online.

The computation of MRPI sets relies on the recursion of predecessor sets, a.k.a. backward reachable sets. Although the number of recursions is theoretically bounded (Rakovic, 2005), it can become computationally challenging, especially for a high dimensional system, to compute the MRPI set at each time step. This computational difficulty arises from the projection operation involved in the predecessor set computation.

For the case where online MRPI set computation is not viable, one can adopt a variant of the proposed online learning based risk-averse MPC that uses set $\mathbb{Z}_f^*$ as the terminal set rather than updating the MRPI set online. According to **Remark 4**, set $\mathbb{Z}_f^*$ can be computed offline and is

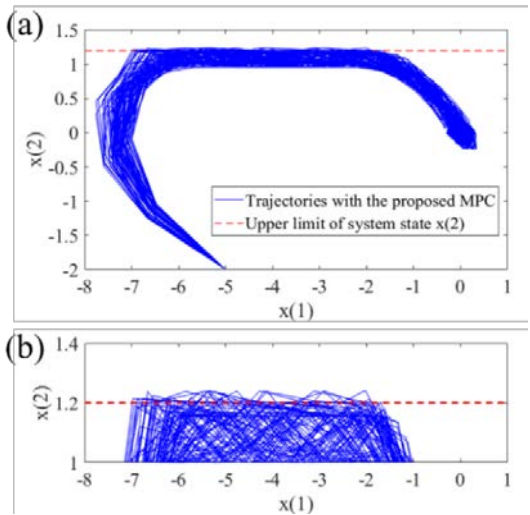

**Fig. 3.** (a): The closed-loop trajectories of system states for the proposed online learning based risk-averse stochastic MPC with 100 realizations of disturbance sequences, (b): The zoomed-in view of state trajectories near the upper limit of $x(2)$.

By leveraging the online updating of the ambiguity set, the time-varying distribution is well captured by the proposed MPC method, and the constraint tightening is adaptive to the distribution accordingly. **Fig. 4** shows the value of parameter $\eta$ in the constraint tightening of the proposed MPC in a simulation run. From **Fig. 4**, we can see that the effect of adaptation in the proposed MPC is evident. Specifically, $\eta$ increases from 0.016 to 0.10 based on the updated information of the stochastic disturbance. The values of *flag* equal one over the entire horizon, which indicates that the adaptive solution is activated at each time step. In this example, the average computational times for solving the finite-horizon optimal control problem, the constraint tightening through (20), and computing an MRPI set at each sampling instant are 0.479s, 0.294s, and 0.123s, respectively. Based on the results, the proposed MPC approach demonstrates its computational efficiency.

always robust positively invariant for the updated tightened constraints. As a result, this modified MPC approach does not require the online computation of MRPI sets, thus saving computational time. In Example 1, it takes 0.057s on average to compute the MRPI set online, and the longest time to compute the MRPI set is 0.105s. By constrast, the modified version of the proposed MPC requires 0.067s offline to compute set $\mathbb{Z}_f^*$, and obviates the computational burden associated with updating MRPI sets online. Additionally, this variant achieves the same average closed-loop cost over all simulation runs compared with the online learning based risk-averse stochastic MPC with updated terminal set $\mathbb{Z}_f^{(k+1)}$. In Example 2, the variant of the proposed MPC method has the same average constraint violation of 7.2% as the original version, and requires only 0.15s offline to compute its terminal set. For general cases, this variant could introduce certain levels of conservatism due to its less flexible construction of the terminal set. Nevertheless, its feature regarding the offline computation of the terminal set is more suitable for high dimensional systems, which is illustrated using the following numerical example (Saltık et al., 2018).

$$x_{k+1} = \begin{pmatrix} 1 & 0 & 0.1 & 0 \\ 0 & 1 & 0 & 0.1 \\ -2 & 0.2 & 1 & 0 \\ 0.5 & -0.05 & 0 & 1 \end{pmatrix} x_k + \begin{pmatrix} 0 \\ 0 \\ 0.2 \\ 0 \end{pmatrix} u_k + w_k \quad (38)$$

The matrix $Q = 5 \cdot I_4$, $R=1$, support set of disturbance $\mathbb{W} = \{w \mid \|w\|_1 \leq 0.06\}$, $\mathbb{U} = \{u \in \mathbb{R} \mid |u| \leq 5\}$, and the distributionally robust CVaR state constraint is $\sup_{\mathbb{P} \in \mathcal{D}^{(k)}} \mathbb{P}\text{-CVaR}_{0.15}\left([0\ 0\ 1\ 0] x_{l+1|k} - 10\right) \leq 0$. The prediction horizon is set to be six, and the simulation horizon is 20. In this example, the online learning based risk-averse stochastic MPC with terminal set $\mathbb{Z}_f^{(k+1)}$ takes an average of 5.511s at each time step to compute the MRPI set in a simulation, which can be computationally expensive for control problems. The modified version needs 5.889s offline to compute the terminal set $\mathbb{Z}_f^*$, and has the same control performance in close-loop cost as the MPC with terminal set $\mathbb{Z}_f^{(k+1)}$.

## 6. Conclusions

In this paper, an online learning based risk-averse stochastic MPC framework was developed for the control of linear time-invariant systems affected by additive disturbance. It incorporated online learning into MPC with desirable theoretical control guarantees and less conservative control performance. Based on the DPMM, a systematic approach to construct the ambiguity set from real-time disturbance data was developed, which leveraged the structural property of multimodality and local moment information. Additionally, the exact reformulation for the distributionally robust CVaR constraint was derived as LMI constraints to facilitate the constraint tightening. The online adaptation of the ambiguity set with real-time measurements helped to improve control performance by actively learning disturbance distributions online. We introduced the safe update scheme to ensure the recursive feasibility and stability of the resulting MPC. The computational results demonstrated that the control performance of the developed MPC is less conservative compared with the one using global mean and covariance information of the disturbance distribution. Additionally, for the case with a time-varying disturbance distribution, the deterioration of constraint satisfaction was ameliorated in the proposed online learning based MPC framework.

## Acknowledgements

The authors acknowledge financial support from the National Science Foundation (NSF) CAREER Award (CBET-1643244). The authors would like to thank the editor, associate editor, and all the reviewers for their valuable and constructive comments, which help improve this paper. We also would like to thank Max Greenberg for his kind help in proofreading the paper.

## Appendix A. The derivation of control objectives

Based on dual-mode prediction dynamics, the nominal system $z_{l+1|k} = Az_{l|k} + Bv_{l|k}$, $z_{0|k} = x_k$ can be characterized by an autonomous system shown below (Kouvaritakis & Cannon, 2016).

$$y_{l+1|k} = \Psi y_{l|k}, \quad l \geq 0 \quad (A.1)$$

where the initial state $y_{0|k} = \begin{bmatrix} z_{0|k} \\ c_{0|k} \\ \vdots \\ c_{N-1|k} \end{bmatrix}$, and the system matrix

$\Psi = \begin{pmatrix} \Phi & B\Xi \\ 0 & \Gamma \end{pmatrix}$, matrix $\Xi = [I_m\ 0 \cdots 0]$, and matrix

$\Gamma = \begin{pmatrix} 0 & I_m & \cdots & 0 \\ \vdots & \vdots & \ddots & \vdots \\ 0 & 0 & \cdots & I_m \\ 0 & 0 & \cdots & 0 \end{pmatrix}$.

Based on the autonomous system dynamics, the stage cost of $J_\infty = \sum_{l=0}^{\infty} \left( z_{l|k}^T Q z_{l|k} + v_{l|k}^T R v_{l|k} \right)$ can be written below.

$$\begin{aligned} z_{l|k}^T Q z_{l|k} + v_{l|k}^T R v_{l|k} &= z_{l|k}^T Q z_{l|k} + \left(K z_{l|k} + c_{l|k}\right)^T R \left(K z_{l|k} + c_{l|k}\right) \\ &= y_{l|k}^T \hat{Q} y_{l|k} \end{aligned} \quad (A.2)$$

where matrix $\hat{Q} = \begin{pmatrix} Q + K^T RK & K^T R\Xi \\ \Xi^T RK & \Xi^T R\Xi \end{pmatrix}$.

Thus, we can rewrite $J_\infty$ as follows.

$$J_\infty = \sum_{l=0}^{\infty} y_{l|k}^T \hat{Q} y_{l|k} = y_{0|k}^T \Theta y_{0|k} \quad (A.3)$$

where $\Theta$ is the positive definite solution of the Lyapunov equation $\Theta = \Psi^T \Theta \Psi + \hat{Q}$.

We express matrix $\Theta$ as $\Theta = \begin{pmatrix} \Theta_z & \Theta_{zc} \\ \Theta_{cz} & \Theta_c \end{pmatrix}$. Substituting $\Theta$, $\Psi$ and $\hat{Q}$ into the equation $\Theta = \Psi^T \Theta \Psi + \hat{Q}$ leads to the following equations:

$$\Theta_z = \Phi^T \Theta_z \Phi + Q + K^T RK \quad (A.4)$$

$$\Theta_{cz} = \Gamma^T \Theta_{cz} \Phi + \Xi^T \left(B^T \Theta_z \Phi + RK\right) \quad (A.5)$$

$$\Theta_c = (B\Xi)^T \Theta_z (B\Xi) + \Gamma^T \Theta_{cz} B\Xi + (B\Xi)^T \Theta_{zc} \Gamma + \Gamma^T \Theta_c \Gamma + \Xi^T R\Xi \tag{A.6}$$

As the feedback gain $K$ is the LQR optimal solution, it follows from (A.4) that $\Theta_z$ is the unique solution of the Riccati equation and $\Theta_z = P$. Based on the expression of LQR solution $K = -(B^T \Theta_z B + R)^{-1} B^T \Theta_z A$, we can further have

$$\begin{aligned}
& B^T \Theta_z \Phi + RK \\
&= B^T \Theta_z \left(A - B(B^T \Theta_z B + R)^{-1} B^T \Theta_z A\right) - R(B^T \Theta_z B + R)^{-1} B^T \Theta_z A \\
&= B^T \Theta_z A - (B^T \Theta_z B + R)(B^T \Theta_z B + R)^{-1} B^T \Theta_z A \\
&= 0
\end{aligned} \tag{A.7}$$

According to (A.5) and (A.7), we have $\Theta_{cz} - \Gamma^T \Theta_{cz} \Phi = 0$, which implies $\Theta_{cz} = 0$ is the solution (Kouvaritakis & Cannon, 2016). As matrix $\Theta$ is symmetric, we can further have $\Theta_{zc} = 0$. Therefore, the equation (A.5) leads to the following equation.

$$\Theta_c - \Gamma^T \Theta_c \Gamma = \Xi^T (R + B^T PB) \Xi \tag{A.8}$$

Let matrix $\Theta_c = \begin{pmatrix} \Theta_c^{(1,1)} & \cdots & \Theta_c^{(1,N)} \\ \vdots & \ddots & \vdots \\ \Theta_c^{(N,1)} & \cdots & \Theta_c^{(N,N)} \end{pmatrix}$. By plugging the block matrix into (A.8), we arrive at the following expression:

$$\Theta_c^{(i,j)} = \begin{cases} R + B^T PB, & i = j \\ 0, & i \neq j \end{cases} \tag{A.9}$$

Based on (A.3) and (A.9), we have $J_\infty = \sum_{l=0}^{N-1} c_{l|k}^T (R + B^T PB) c_{l|k} + x_k^T P x_k$.

**Appendix B. Proof of Theorem 1**

**Proof.** Based on $x_{l|k} = z_{l|k} + e_{l|k}$ and $e_{l+1|k} = \Phi e_{l|k} + w_{l|k}$, we have

$$x_{l+1|k} = z_{l+1|k} + \Phi e_{l|k} + w_{l|k} \tag{B.1}$$

Note that conditional distributionally robust CVaR constraints need to hold for any reachable $e_{l|k}$. Therefore, we leverage the tube-based method, and the original distributionally robust CVaR constraints can be cast as follows.

$$z_{l+1|k} \in \hat{\mathbb{Z}}_{l+1}^{(k)} = \tilde{\mathbb{Z}}^{(k)} \ominus \left(\bigoplus_{i=1}^{l} \Phi^i \mathbb{W}\right) \tag{B.2}$$

where set $\tilde{\mathbb{Z}}^{(k)}$ is explicitly given by,

$$\tilde{\mathbb{Z}}^{(k)} = \left\{ z \middle| \sup_{\mathbb{P} \in \mathcal{D}^{(k)}} \mathbb{P}\text{-CVaR}_{[\varepsilon]_i} \left([H]_i z - [h]_i + [H]_i w_{l|k}\right) \leq 0, i \in [1, p] \right\} \tag{B.3}$$

We will provide a constraint tightening method to reformulate the following constraint.

$$\sup_{\mathbb{P} \in \mathcal{D}^{(k)}} \mathbb{P}\text{-CVaR}_{[\varepsilon]_i} \left([H]_i z - [h]_i + [H]_i w_{l|k}\right) \leq 0 \tag{B.4}$$

We can further reformulate the distributionally robust CVaR constraint, according to the definition of CVaR and a stochastic minimax theorem (Hanasusanto et al., 2015; Shapiro & Kleywegt, 2002; Zymler et al., 2013).

$$\begin{aligned}
& \sup_{\mathbb{P} \in \mathcal{D}^{(k)}} \mathbb{P}\text{-CVaR}_{[\varepsilon]_i} \left([H]_i z - [h]_i + [H]_i w_{l|k}\right) \\
&= \sup_{\mathbb{P} \in \mathcal{D}^{(k)}} \inf_{\beta_i \in \mathbb{R}} \left\{\beta_i + \frac{1}{[\varepsilon]_i} \mathbb{E}_\mathbb{P}\left\{\left([H]_i z - [h]_i + [H]_i w_{l|k} - \beta_i\right)^+\right\}\right\} \\
&= \inf_{\beta_i \in \mathbb{R}} \left\{\beta_i + \frac{1}{[\varepsilon]_i} \sup_{\mathbb{P} \in \mathcal{D}^{(k)}} \mathbb{E}_\mathbb{P}\left\{\left([H]_i z - [h]_i + [H]_i w_{l|k} - \beta_i\right)^+\right\}\right\}
\end{aligned} \tag{B.5}$$

where $\mathbb{E}_\mathbb{P}\{\cdot\}$ represents the expectation with respect to probability $\mathbb{P}$. The first equality holds based on **Definition 1**. The second equality is based on the stochastic minimax theorem.

The worst-case expectation problem $\sup_{\mathbb{P} \in \mathcal{D}^{(k)}} \mathbb{E}_\mathbb{P}\left\{\left([H]_i z - [h]_i + [H]_i w_{l|k} - \beta_i\right)^+\right\}$ can be rewritten as the following optimization problem.

$$\begin{aligned}
& \sup_{\rho_1, \ldots, \rho_m \in \mathcal{M}_+} \sum_{j=1}^{m} \gamma_j \int_\mathbb{W} \left([H]_i z - [h]_i + [H]_i \xi - \beta_i\right)^+ \rho_j(d\xi) \\
& \text{s.t.} \quad \begin{cases} \int_\mathbb{W} \rho_j(d\xi) = 1 \\ \int_\mathbb{W} \xi \, \rho_j(d\xi) = \mu_j \\ \int_\mathbb{W} \xi \xi^T \rho_j(d\xi) \leq \Sigma_j + \mu_j \mu_j^T \end{cases} \forall j = 1, \ldots, m
\end{aligned} \tag{B.6}$$

By taking the dual of optimization problem (B.6), we have (B.7)-(B.10).

$$\min_{t_{ij}, \omega_{ij}, \Omega_{ij}} \sum_{j=1}^{m} \gamma_j \left\{t_{ij} + \mu_j^T \omega_{ij} + \left(\Sigma_j + \mu_j \mu_j^T\right) \bullet \Omega_{ij}\right\} \tag{B.7}$$

$$\text{s.t.} \quad t_{ij} + \xi^T \omega_{ij} + \xi^T \Omega_{ij} \xi \geq 0, \ \forall \xi \in \mathbb{W}, \ \forall j \tag{B.8}$$

$$\begin{aligned}
& -\left([H]_i z - [h]_i + [H]_i \xi - \beta_i\right) + t_{ij} \\
& + \xi^T \omega_{ij} + \xi^T \Omega_{ij} \xi \geq 0, \ \forall \xi \in \mathbb{W}, \ \forall j
\end{aligned} \tag{B.9}$$

$$\Omega_{ij} \geq 0, \ \forall j \tag{B.10}$$

where $t_{ij}$, $\omega_{ij}$, and $\Omega_{ij}$ are the dual variables corresponding to the constraints in the ambiguity set. Constraints (B.8)-(B.9) are semi-infinite constraints, so some further reformulations are needed.

Constraint $t_{ij} + \xi^T \omega_{ij} + \xi^T \Omega_{ij} \xi \geq 0, \ \forall \xi \in \mathbb{W}, \ \forall j$ can be reformulated as follows.

$$\min_{\xi \in \mathbb{W}} \left\{t_{ij} + \xi^T \omega_{ij} + \xi^T \Omega_{ij} \xi\right\} \geq 0, \ \forall j \tag{B.11}$$

We can further reformulate constraint $\min_{\xi \in \mathbb{W}}\{t_{ij} + \xi^T \omega_{ij} + \xi^T \Omega_{ij} \xi\} \geq 0, \forall j$ as the following LMI constraints, based on the duality of convex quadratic programs and Schur complements (Boyd & Vandenberghe, 2004).

$$\begin{pmatrix} \Omega_{ij} & \frac{1}{2}(\omega_{ij} + E^T \varphi_{ij}) \\ \frac{1}{2}(\omega_{ij} + E^T \varphi_{ij})^T & t_{ij} - f^T \varphi_{ij} \end{pmatrix} \geq 0, \forall j=1,\ldots,m \quad (B.12)$$

where $\varphi_{ij} \geq 0$ is the vector of dual variables corresponding to the constraints $E\xi \leq f$.

Similarly, constraints (B.9) can be recast as LMI constraints below.

$$\begin{pmatrix} \Omega_{ij} & \frac{1}{2}(\omega_{ij} - [H]_i^T + E^T \phi_{ij}) \\ \frac{1}{2}(\omega_{ij} - [H]_i^T + E^T \phi_{ij})^T & t_{ij} + \beta_i + [h]_i - [H]_i z - f^T \phi_{ij} \end{pmatrix} \geq 0, \forall j \quad (B.13)$$

where $\phi_{ij} \geq 0$ is the vector of dual variables corresponding to the constraints $E\xi \leq f$.

Thus, the distributionally robust CVaR constraints $\sup_{\mathbb{P} \in \mathcal{D}^{(k)}} \mathbb{P}\text{-CVaR}_{[\varepsilon]_i}([H]_i z - [h]_i + [H]_i w_{l|k}) \leq 0$ over the DPMM-based ambiguity set are reformulated as follows.

$$[\varepsilon]_i \cdot \beta_i + \sum_{j=1}^m \gamma_j \{t_{ij} + \mu_j^T \omega_{ij} + (\Sigma_j + \mu_j \mu_j^T) \bullet \Omega_{ij}\} \leq 0$$

$$\begin{pmatrix} \Omega_{ij} & \frac{1}{2}(\omega_{ij} + E^T \varphi_{ij}) \\ \frac{1}{2}(\omega_{ij} + E^T \varphi_{ij})^T & t_{ij} - f^T \varphi_{ij} \end{pmatrix} \geq 0, \forall j \quad (B.14)$$

$$\begin{pmatrix} \Omega_{ij} & \frac{1}{2}(\omega_{ij} - [H]_i^T + E^T \phi_{ij}) \\ \frac{1}{2}(\omega_{ij} - [H]_i^T + E^T \phi_{ij})^T & t_{ij} + \beta_i + [h]_i - [H]_i z - f^T \phi_{ij} \end{pmatrix} \geq 0, \forall j$$

$$\Omega_{ij} \geq 0, \varphi_{ij} \geq 0, \phi_{ij} \geq 0, \forall j$$

We consider the following constraint:

$$[H]_i z \leq [h]_i - [\eta^{(k)}]_i \quad (B.15)$$

where

$$[\eta^{(k)}]_i = \min_\eta \eta$$
$$\text{s.t.} \sup_{\mathbb{P} \in \mathcal{D}^{(k)}} \mathbb{P}\text{-CVaR}_{[\varepsilon]_i}([H]_i w_{l|k} - \eta) \leq 0 \quad (B.16)$$

The constraint in (B.16) can be converted into the constraints in (20). Then, we convert (B.15) and (B.16) as follows.

$$[H]_i z \leq \max_{[\eta^{(k)}]_i \in \{\eta | \sup_{\mathbb{P} \in \mathcal{D}^{(k)}} \mathbb{P}\text{-CVaR}_{[\varepsilon]_i}([H]_i w_{l|k} - \eta) \leq 0\}} ([h]_i - [\eta^{(k)}]_i) \quad (B.17)$$

We reformulate constraint (B.17) into (B.14) using the same reformulation technique. This completes the proof. □

## Appendix C. Proof of Proposition 1

**Proof.** First, we introduce an ambiguity set $\mathcal{D}^* = \{\rho \in \mathcal{M}_+ \mid \int_\mathbb{W} \rho(d\xi) = 1\}$. Consider the constraint tightening based on the ambiguity set $\mathcal{D}^*$ below.

$$[\tilde{\eta}]_i = \min \eta$$
$$\text{s.t.} \sup_{\mathbb{P} \in \mathcal{D}^*} \mathbb{P}\text{-CVaR}_{[\varepsilon]_i}([H]_i w_{l|k} - \eta) \leq 0 \quad (C.1)$$

By using a similar approach as presented in Appendix B, we can reformulate (C.1) to the following problem.

$$[\tilde{\eta}]_i = \min \eta$$
$$\text{s.t.} [\varepsilon]_i \cdot \beta_i + t_i \leq 0$$
$$t_i \geq 0 \quad (C.2)$$
$$t_i + \eta + \beta_i \geq \max_{\xi \in \mathbb{W}} [H]_i \xi$$

According to constraints in (C.2), we can further have $\eta \geq \max_{\xi \in \mathbb{W}} [H]_i \xi + \left(\frac{1}{[\varepsilon]_i} - 1\right) t_i$. Because $t_i \geq 0$, we have $[\tilde{\eta}]_i = \max_{\xi \in \mathbb{W}} [H]_i \xi = [\eta_0]_i$.

Because $\mathcal{D}^{(k)} \subseteq \mathcal{D}^*$, it holds that $[\eta_0]_i \geq [\eta^{(k)}]_i, \forall i, k$. Given inequality $[\eta_0]_i \geq [\eta^{(k)}]_i$, we can further have $\hat{Z}_l^* \subseteq Z_l^{(k)}, \forall k$. Based on $\hat{Z}_l^* \subseteq Z_l^{(k)}, \forall k$ and **Definition 4**, $Z_f^*$ is a robust positively invariant set (not necessarily MRPI set) for the updated tightened constraints. As $Z_f^*$ is nonempty under **Assumption 5**, there always exists a nonempty MRPI set according to **Definition 5**. This completes the proof. □

## Appendix D. Proof of Theorem 2

**Proof.** Given that $\mathbb{D}(x_k) \neq \emptyset$, we want to show $\mathbb{D}(x_{k+1}) \neq \emptyset$. Suppose $\mathbf{c}_{N|k}^* = (c_{0|k}^*, c_{1|k}^*, \ldots, c_{N-1|k}^*)$ is the optimal solution to problem (**OL-SMPC$_k$**) at time $k$. Hence, the corresponding optimal nominal state evolution is given below.

$$z_{l+1|k}^* = \Phi z_{l|k}^* + B c_{l|k}^*, \quad z_{0|k}^* = x_k \quad (D.1)$$

Based on the optimal solution, an explicit candidate solution can be constructed as $\tilde{\mathbf{c}}_{k+1} = (c_{1|k}^*, \ldots, c_{N-1|k}^*, 0)$ using **Definition 3**. The nominal state under control input $\tilde{\mathbf{c}}_{k+1}$ is shown as follows.

$$\tilde{z}_{0|k+1} = x_{k+1} = Ax_k + Bu_k + w_k$$
$$= (A + BK)x_k + Bc_{0|k}^* + w_k \quad (D.2)$$
$$= z_{1|k}^* + w_k$$

Based on $\tilde{z}_{0|k+1} = z^*_{1|k} + w_k$, we further have the following equalities.

$$\begin{aligned}\tilde{z}_{1|k+1} &= \Phi \tilde{z}_{0|k+1} + B\tilde{c}_{0|k+1}\\ &= \Phi\left(z^*_{1|k} + w_k\right) + Bc^*_{1|k}\\ &= z^*_{2|k} + \Phi w_k\end{aligned} \quad (D.3)$$

According to $\tilde{z}_{0|k+1} = z^*_{1|k} + w_k$ and $\tilde{z}_{1|k+1} = z^*_{2|k} + \Phi w_k$, we have the following relation by induction.

$$\tilde{z}_{l|k+1} = z^*_{l+1|k} + \Phi^l w_k, \quad \forall l \in [0, N-1] \quad (D.4)$$

Similarly, we have the relationship between the optimal solution at time $k$ and the candidate at time $k+1$, as presented by,

$$\tilde{v}_{l|k+1} = v^*_{l+1|k} + K\Phi^l w_k, \quad \forall l \in [0, N-1] \quad (D.5)$$

For the optimal control problem (**OL-SMPC$_{k+1}$**), we can consider two scenarios, namely *flag*=0 and *flag*=1. For the scenario in which *flag*=0, we have $\mathbb{Z}_l^{(k+1)} = \mathbb{Z}_l^{(k)}$ according to $\mathbb{Z}_l^{(k+1)} = flag \cdot \hat{\mathbb{Z}}_l^{(k+1)} \oplus (1-flag) \cdot \mathbb{Z}_l^{(k)}$. Based on $\tilde{z}_{l|k+1} = z^*_{l+1|k} + \Phi^l w_k$ and $\tilde{v}_{l|k+1} = v^*_{l+1|k} + K\Phi^l w_k$, we have the following:

$$z^*_{l+1|k} \in \mathbb{Z}^{(k)}_{l+1} \Rightarrow \tilde{z}_{l|k+1} \in \mathbb{Z}^{(k+1)}_l \quad (D.6)$$

$$v^*_{l+1|k} \in \mathbb{V}_{l+1} \Rightarrow \tilde{v}_{l|k+1} \in \mathbb{V}_l \quad (D.7)$$

Next, we derive the candidate nominal state at the end of the horizon as follows.

$$\begin{aligned}\tilde{z}_{N|k+1} &= A\tilde{z}_{N-1|k+1} + B\tilde{v}_{N-1|k+1}\\ &= A\left(z^*_{N|k} + \Phi^{N-1} w_k\right) + B\left(Kz^*_{N|k} + K\Phi^{N-1} w_k\right)\\ &= \Phi z^*_{N|k} + \Phi^N w_k\end{aligned} \quad (D.8)$$

Note that $\mathbb{Z}_f^{(k+1)} = \mathbb{Z}_f^{(k)}$ because *flag*=0. Based on the definition of the terminal set, we have the following.

$$z^*_{N|k} \in \mathbb{Z}^{(k)}_f \Rightarrow \tilde{z}_{N|k+1} \in \mathbb{Z}^{(k+1)}_f \quad (D.9)$$

Up to now, we have checked all the constraints, we can conclude that the candidate solution is feasible for the problem (**OL-SMPC$_{k+1}$**) if *flag*=0.

For the scenario where *flag*=1, it implies that the constructed solution satisfies the constraints in (**OL-SMPC$_{k+1}$**), according to the safe update scheme in Section 3.3, which completes the proof. □

**Appendix E. Proof of Theorem 3**

**Proof.** We define the optimal objective value for the problem (**OL-SMPC$_k$**) as $J_k$ below.

$$J_k \doteq V_N^*(x_k) = \sum_{l=0}^{N-1} \left\|c^*_{l|k}\right\|^2_{\tilde{\Psi}=R+B^T PB} \quad (E.1)$$

Since the candidate solution is feasible and the objective function remains the same for both scenarios where *flag*=0 and *flag*=1, we have

$$J_{k+1} \leq \tilde{V}_N(x_{k+1}) = \sum_{l=0}^{N-1}\left\|\tilde{c}_{l|k+1}\right\|^2_{\tilde{\Psi}} = \sum_{l=1}^{N-1}\left\|c^*_{l|k}\right\|^2_{\tilde{\Psi}} = J_k - \left\|c^*_{0|k}\right\|^2_{\tilde{\Psi}} \quad (E.2)$$

where $\tilde{V}_N(x_{k+1})$ represents the objective value corresponding to the candidate solution.

By rearranging $J_{k+1} \leq J_k - \left\|c^*_{0|k}\right\|^2_{\tilde{\Psi}}$, we have the following inequality.

$$J_{k+1} - J_k \leq -\left\|c^*_{0|k}\right\|^2_{\tilde{\Psi}} \quad (E.3)$$

By adding $J_{k+1} - J_k \leq -\left\|c^*_{0|k}\right\|^2_{\tilde{\Psi}}$ from $k$=0, we have the follows.

$$\sum_{k=0}^{\infty}\left\|c^*_{0|k}\right\|^2_{\tilde{\Psi}} \leq J_0 - J_\infty \quad (E.4)$$

Based on $\sum_{k=0}^{\infty}\left\|c^*_{0|k}\right\|^2_{\tilde{\Psi}} \leq J_0 - J_\infty$, we have $\lim_{k \to \infty}\left\|c^*_{0|k}\right\| = 0$ for both scenarios where *flag*=0 and *flag*=1. Additionally, the convergence of the state to a neighborhood of the origin is further established.

Then, we consider the asymptotical dynamic behavior of the system state as follows.

$$\begin{aligned}\lim_{k \to \infty} x_k &= \lim_{k \to \infty}\left(\Phi^k x_0 + \sum_{i=1}^{k}\Phi^{i-1} B c^*_{0|k-i} + \sum_{i=1}^{k}\Phi^{i-1} w_{k-i}\right)\\ &= \lim_{k \to \infty}\left(\sum_{i=1}^{k}\Phi^{i-1} w_{k-i}\right)\end{aligned} \quad (E.5)$$

Note that the second equality holds because $\Phi$ is Shur and $\lim_{k \to \infty}\left\|c^*_{0|k}\right\| = 0$.

According to $\lim_{k \to \infty} x_k = \lim_{k \to \infty}\left(\sum_{i=1}^{k}\Phi^{i-1} w_{k-i}\right)$, we have the following relation.

$$\lim_{k \to \infty} x_k \in R_\infty \doteq \bigoplus_{i=0}^{\infty} \Phi^i \mathbb{W} \quad (E.6)$$

Based on $\lim_{k \to \infty} x_k \in R_\infty$, the system state converges to a neighborhood of the origin, namely the minimal robust positively invariant set, under the proposed online learning based risk-averse stochastic MPC strategy. This completes the proof. □